\documentclass[reqno,12pt,fleqn]{amsart}
\usepackage{amsmath,amsthm,amsfonts,amssymb,latexsym,epsfig}
\usepackage{times,csquotes}

\usepackage{geometry,float,setspace}
\usepackage{cite}
\usepackage[colorlinks=true, allcolors=blue]{hyperref}
\usepackage{cleveref}
\crefname{equation}{}{}
\crefformat{equation}{(#2#1#3)}
\crefrangeformat{equation}{(#3#1#4)--(#5#2#6)}
\allowdisplaybreaks
\Crefformat{equation}{Equation~(#2#1#3)}
\newtheorem{theorem}{Theorem}[section]

\newtheorem{definition}{Definition}[section]
\newtheorem{lemma}{Lemma}[section]

\theoremstyle{remark}
\theoremstyle{definition}

\numberwithin{equation}{section}

\newcommand\numberthis{\addtocounter{equation}{1}\tag{\theequation}}
\makeatletter
\def\specialsection{\@startsection{section}{1}
	\z@{\linespacing\@plus\linespacing}{.5\linespacing}
	{\normalfont}}
\def\section{\@startsection{section}{1}
	\z@{.7\linespacing\@plus\linespacing}{.5\linespacing}
	{\normalfont\scshape}}
\everymath{\displaystyle{}}
\usepackage{xparse}
\NewDocumentCommand{\p}{r()}{
	\ensuremath{\left(#1\right)_{\infty}}
}
\newenvironment{nalign}{%
  \align
}{%
  \endalign
}
\newenvironment{nalign*}{%
  \csname align*\endcsname
  \setlength{\mathindent}{0pt}%
}{%
  \endalign
}
\AtBeginEnvironment{nalign}{\setlength{\mathindent}{0pt}}
\AtBeginEnvironment{nalign*}{\setlength{\mathindent}{0pt}}
\title{R3R}
\begin{document}
\title[{\tiny  Vanishing coefficient results in four families of infinite $q$-products }]{Vanishing coefficient results in four families of infinite $q$-products  }
 \maketitle
 \begin{center}{S. Ananya$^1$, Channabasavayya$^1$, D. Ranganatha$^\star$ and R. G. Veeresha$^2$  \\
 \vspace{.5cm}
  $^{1, \star}$Department of Mathematics, Central University of Karnataka, India\\
  $^2$Department of Mathematics, Sri Jayachamarajendra College of Engineering,\\ JSS Science and Technology University, Mysuru- 570006, Karnataka, India\\
  \vspace{.5cm}
E-mail: ananyasahadev2525@gmail.com, cshmatics@gmail.com, ddranganatha@gmail.com$^\star$\\  and  veeru.rg@gmail.com}
 \end{center}
\begin{abstract}
    In the recent past, the work in the area of vanishing coefficients of infinite $q$-products has been taken to the forefront. Weaving the same thread as  Ramanujan, Richmond, Szekeres, Andrews, Alladi, Gordon, Mc Laughlin, Baruah, Kaur, Tang,  we further prove vanishing coefficients in arithmetic progressions  moduli 5, 7, 11, 13, 19, 21, 23 and 29 of the following four families of infinite products, where $\{X_{a,b,sm,km,u,v}(n)\}_{n\geq n_0}$, $\{Y_{a,b,sm,km,u,v}(n)\}_{n\geq n_0}$, $\{Z_{a,b,sm,km,u,v}(n)\}_{n\geq n_0}$ and $\{W_{a,b,sm,km,u,v}(n)  \}_{n\geq n_0}$ are defined by
\begin{align*}
\sum_{n\geq n_0}^{\infty}X_{a,b,sm,km,u,v}(n)q^n:=&\p(q^{a},q^{sm-a};q^{sm})^u\p(q^{b},q^{km-b};q^{km})^v, \\
\sum_{n\geq n_0}^{\infty}Y_{a,b,sm,km,u,v}(n)q^n:=&\p(q^{a},q^{sm-a};q^{sm})^u\p(-q^{b},-q^{km-b};q^{km})^v, \\
\sum_{n\geq n_0}^{\infty}Z_{a,b,sm,km,u,v}(n)q^n:=&\p(-q^{a},-q^{sm-a};q^{sm})^u\p(q^{b},q^{km-b};q^{km})^v,\\
  \sum_{n\geq n_0}^{\infty}W_{a,b,sm,km,u,v}(n)q^n:=&\p(-q^{a},-q^{sm-a};q^{sm})^u\p(-q^{b},-q^{km-b};q^{km})^v,
\end{align*}
here $a, b, s, k, u$ and $v$  are chosen in such a way that the infinite products in the right-hand side of the above are convergent and  $n_0$ is an integer (possibly
negative or zero) depending on $a, b, s, k, u$ and $v$. The proof uses the Jacobi triple product identity and the properties of Ramanujan general theta function. \\

\noindent\textsc{Mathematics Subject Classification.} 11B65, 11F27, 30B10.\\

 \noindent\textsc{Keywords and phrases.}
Vanishing coefficients, products of theta functions, arithmetic
progressions, odd moduli.
\end{abstract}
\maketitle
\section{Introduction}
The reference \cite[p. 49]{ramanujan1988lost} provides insights into the calculation of coefficients of the following infinite $q$-products done by Ramanujan for up to $n=1000$, but provides no conclusion:
\begin{align*}
    \sum_{n=0}^{\infty}u(n)q^n:=\frac{ \left(q^4 ; q^5\right)_{\infty}\left(q; q^5\right)_{\infty}}{\left(q^2; q^5\right)_{\infty} \left(q^3 ; q^5\right)_{\infty}} \quad \textrm{and}\quad \sum_{n=0}^{\infty}v(n)q^n:=\left(\frac{ \left(q^3 ; q^5\right)_{\infty}\left(q^2; q^5\right)_{\infty}}{\left(q; q^5\right)_{\infty} \left(q^4 ; q^5\right)_{\infty}}\right)^{-1},
\end{align*} 
here and throughout the paper, we use the following $q$-product notations:
\begin{align*}
    (a;q)_\infty&=\prod_{k=1}^\infty(1-aq^{k-1}) \quad (a, q\in\mathbb{C},\quad |q|<1)
    \end{align*}
    and
    \begin{align*}
(a_1,a_2,\hdots,a_m;q)_\infty&=(a_1;q)_\infty(a_2;q)_\infty\cdots(a_m;q)_\infty, \quad \textrm{where}  \quad a_1, a_2, a_3,\hdots, a_m \in \mathbb{C}.
\end{align*} 
Motivated by the works of Ramanujan, Richmond and Szekeres \cite{richmond1978taylor} explored the vanishing coefficients in the following infinite products and  proved that if
$$
\sum_{n=0}^{\infty} \alpha(n) q^n:=\frac{\left(q^3;q^8\right)_\infty \left(q^5 ; q^8\right)_{\infty}}{\left(q;q^8\right)_\infty \left(q^7 ; q^8\right)_{\infty}} \quad \text { and } \quad \sum_{n=0}^{\infty} \beta(n) q^n:=\frac{\left(q;q^8\right)_\infty \left(q^7 ; q^8\right)_{\infty}}{\left(q^3;q^8\right)_\infty \left(q^5 ; q^8\right)_{\infty}},
$$
then $\alpha(4n+3)=0$ and $\beta(4n+2)=0.$ Along with it, they proposed the following conjecture: 
$$
\sum_{n=0}^{\infty} \gamma(n) q^n:=\frac{\left(q^5;q^{12}\right)_\infty \left(q^7 ; q^{12}\right)_{\infty}}{\left(q;q^{12}\right)_\infty \left(q^{11} ; q^{12}\right)_{\infty}} \quad\text { and } \quad \sum_{n=0}^{\infty} \delta(n) q^n:=\frac{\left(q;q^{12}\right)_\infty \left(q^{11} ; q^{12}\right)_{\infty}}{\left(q^5;q^{12}\right)_\infty \left(q^7 ; q^{12}\right)_{\infty}},
$$
then $\gamma(6 n+5)=0$ and $\delta(6 n+3)=0$.  This led to the work of Andrews and Bressoud \cite{Andrewsvc}. They proved  the following  theorem using Ramanujan's ${}_1\psi_1$ formula, which included Richmond and Szekeres' results \cite{richmond1978taylor} and their conjecture as special cases: 
\begin{theorem}
 If $1 \leq r<k$ are relatively prime integers of opposite parity and
$$
\sum_{n=0}^{\infty} \phi(n) q^n:=\frac{(q^r;q^{2k})_\infty(q^{2 k-r};q^{2 k})_{\infty}}{(q^{k-r};q^{2k})_\infty (q^{k+r};q^{2 k})_{\infty}},
$$
then $\phi(k n+r(k-r+1) / 2)=0$. \end{theorem}
The results of these types were further generalized by Alladi and Gordon \cite{Alladi} and James Mc Laughlin \cite{Mc3}.

\par
In 2018, Hirschhorn \cite{hirschhorn2019two} proved the following theorem, which serves as a foundational step in enriching knowledge on vanishing coefficients in other classes of infinite products:
\begin{theorem}  Let the sequences \(\{a(n)\}\) and \(\{b(n)\}\) be defined by
\begin{align}
\sum_{n=0}^{\infty}a(n)q^n := (-q,-q^4;q^5)_\infty (q, q^9;q^{10})^3_\infty, \label{MDHvc1}\\
\sum_{n=0}^{\infty}b(n)q^n := (-q^2,-q^3;q^5)_\infty (q^3, q^7;q^{10})^3_\infty, \label{MDHvc2}
\end{align} then
$a(5n+2)=a(5n+4)=b(5n+1)=b(5n+4)=0.$   
\end{theorem}
The proof employed $q$-series manipulations, Jacobi's triple product identity (JTPI), and linear transformations to the indices. \par
In the same year, Tang \cite{tang2019vanishing} considered the variants of \eqref{MDHvc1} and \eqref{MDHvc2} obtained some comparable results. For example, he proved that if \(\{a_1(n)\}\), \(\{b_1(n)\}\), \(\{a_2(n)\}\) and \(\{b_2(n)\}\)  are defined as below,
\begin{nalign*}
    \sum_{n=0}^{\infty} a_1(n) q^n &:= (-q,-q^4;q^5)^3_\infty  (q^2,q^8;q^{10})_\infty,\quad
        \sum_{n=0}^{\infty}b_1(n)q^n := (-q^2,-q^3;q^5)^3_\infty  (q^4,q^6;q^{10})_\infty,\\
    \sum_{n=0}^{\infty} a_2(n) q^n &:= (-q,-q^4;q^5)^3_\infty  (q^3,q^7;q^{10})_\infty,\quad
        \sum_{n=0}^{\infty}b_2(n) q^n:= (-q^2,-q^3;q^5)^3_\infty  (q,q^9;q^{10})_\infty,
    \end{nalign*}
then 
        \begin{align*}
        a_1(5n+4) = b_1(5n+1) = a_2(5n+3) = a_2(5n+4) = b_2(5n+3) = b_2(5n+4) = 0.
        \end{align*}
        
 Soon after, Baruah and Kaur \cite{baruah2020some} obtained some new results of the above type and proved them using JTPI and elementary $q$-series manipulation techniques. 

 In the recent past, Mc Laughlin \cite{mc2021new} studied and presented that vanishing coefficients of some infinite products can be categorized into families and they can proved. To cite his work, he proved if $t \in\{1,2\}$ define the sequence $\left\{s(n)\right\}$ by
$$
\sum_{n=0}^{\infty} s(n) q^n:=(q^{2t},q^{5-2t};q^5)_{\infty}(q^{5-2t},q^{5+2t} ; q^{10})_{\infty}^3, 
$$
then $s({5n+3t^2+t})=0$ for all $n$. \par
Mc Laughlin used the Quintuple product identity (QPI) \cite[Eq. (4.6)]{Mc3} to prove vanishing coefficients in the infinite products, whereas others have used JTPI and properties of  Ramanujan's general theta function. He also extended such vanishing coefficient results to other modulus, such as 7 and 11. In \cite{mc2021new}, Mc Laughlin proved several results and listed several results on vanishing coefficients in infinite products with negative signs, which were proved by Kaur and Vandna \cite{vm1, vm2} quite recently. 
In 2023, Tang \cite{2tang2023vanishing}  proved many vanishing coefficient identities of  three infinite products  by using properties of Ramanujan's theta function,  QPI and an extended QPI (EQPI) due to Cao \cite{cao2011integer}, and these types of products were studied by many researchers, as mentioned above. He generalized the results of  Hirschhorn \cite{hirschhorn2019two}, Vandna and Kaur \cite{vm1, vm2}, Baruah and Kaur \cite{baruah2020some}, Mc Laughlin \cite{mc2021new}, etc. He provided a total of 36 results and offered two conjectures. Mc Laughlin \cite{mc2021new} results could be applied to only a finite number of specific cases. But, these new results due to Tang \cite{2tang2023vanishing} can be applied to an infinite number of particular cases. In \cite{tang2023vanishing},  Tang continued the same for odd moduli. Finally, he posed several conjectures that were still not answered. For More recent work, please see \cite{dasappa2024vanishing, Ranganatha, keerthana2024generalization}. 
Motivated by the above work, in this paper we define the following four infinite products:
\begin{definition}
Define $\{X_{a,b,s\ell,k\ell,u,v}(n)\}_{n\geq n_0}$, $\{Y_{a,b,s\ell,k\ell,u,v}(n)\}_{n\geq n_0}$, $\{Z_{a,b,s\ell,k\ell,u,v}(n)\}_{n\geq n_0}$ and \\$\{W_{a,b,s\ell,k\ell,u,v}(n)\}_{n\geq n_0}$ by
\begin{align}
\label{D1}\sum_{n\geq n_0}^{\infty}X_{a,b,s\ell,k\ell,u,v}(n)q^n=&\p(q^{a},q^{s\ell-a};q^{s\ell})^u\p(q^{b},q^{k\ell-b};q^{k\ell})^v, \\
\label{D2}\sum_{n\geq n_0}^{\infty}Y_{a,b,s\ell,k\ell,u,v}(n)q^n=&\p(q^{a},q^{s\ell-a};q^{s\ell})^u\p(-q^{b},-q^{k\ell-b};q^{k\ell})^v, \\
\label{D3}\sum_{n\geq n_0}^{\infty}Z_{a,b,s\ell,k\ell,u,v}(n)q^n=&\p(-q^{a},-q^{s\ell-a};q^{s\ell})^u\p(q^{b},q^{k\ell-b};q^{k\ell})^v,
\end{align}
and
\begin{align}
  \label{D4}  \sum_{n\geq n_0}^{\infty}W_{a,b,s\ell,k\ell,u,v}(n)q^n=&\p(-q^{a},-q^{s\ell-a};q^{s\ell})^u\p(-q^{b},-q^{k\ell-b};q^{k\ell})^v,
\end{align}
where $a, b, s, k, \ell,  u$ and $v$  are chosen in such a way that the infinite products in right-hand side of (\ref{D1})--(\ref{D4}) are convergent and  $n_0$ is an integer (possibly
negative or zero) depending on $a, b, s, k, u$ and $v$.
\end{definition}
Historically, the coefficient-vanishing property for the above infinite products with $s=1, k=2$ and a few values of $u, v$ has been studied in the above. In this paper, we investigate the vanishing coefficients for $s=k=1$,  $s=1, k=3$ and some values of $u$ and $v$. In order to state our results, we adapt the following convention here.  \begin{align}
\label{C1}X_{a,b,s\ell,k\ell,u,v}(pn+qt)=0~ \textrm{holds ~ under~ the ~ condition ~ } ~\gcd(p, t)=1,
\end{align}
where $p$, $q$ are some positive integers and  $t$ is any integer.

The main results of this paper are stated as follows and we use the JTPI and the properties of Ramanujan general theta function to establish the same results:
\begin{theorem}\label{vcthm2} Let $\ell \geq 1$ and $t$ satisfying the condition stated in (\ref{C1}). Then 
\begin{align}
&X_{t,2t,5\ell,15\ell,2,1}(5n+2t)=Z_{t,2t,5\ell,15\ell,2,1}(5n+2t)=0,\label{vcres2.0}\\ 
&X_{2t,t,7\ell,21\ell,1,2}(7n+2t)=Y_{2t,t,7\ell,21\ell,1,2}(7n+2t)=0,\label{vcres2.1}\\ 
&X_{t,6t,7\ell,21\ell,2,1}(7n+4t)=Z_{t,6t,7\ell,21\ell,2,1}(7n+4t)=0,\label{vcres2.2}\\ 
&X_{4t,3t,11\ell,33\ell,1,2}(11n+5t)=Y_{4t,3t,11\ell,33\ell,1,2}(11n+5t)=0,\label{vcres2.3}\\ 
&X_{2t,8t,11\ell,33\ell,2,1}(11n+6t)=Z_{2t,8t,11\ell,33\ell,2,1}(11n+6t)=0,\label{vcres2.4}\\ 
&X_{2t,3t,13\ell,39\ell,1,3}(13n+12t)=Y_{2t,3t,13\ell,39\ell,1,3}(13n+12t)=0,\label{vcres2.9}\\
&X_{2t,9t,13\ell,39\ell,3,1}(13n+t)=Y_{2t,9t,13\ell,39\ell,3,1}(13n+t)=0,\label{vcres2.10}\\ 
&X_{t,12t,13\ell,39\ell,4,1}(13n+8t)=Z_{t,12t,13\ell,39\ell,4,1}(13n+8t)=0,\label{vcres2.11}\\ 
&X_{4t,t,13\ell,39\ell,1,4}(13n+4t)=Y_{4t,t,13\ell,39\ell,1,4}(13n+4t)=0,\label{vcres2.12}\\ 
&X_{8t,3t,19\ell,57\ell,1,4}(19n+10t)=Y_{8t,3t,19\ell,57\ell,1,4}(19n+10t)=0,\label{vcres2.16}\\
&X_{t,8t,19\ell,57\ell,4,1}(19n+6t)=Z_{t,8t,19\ell,57\ell,4,1}(19n+6t)=0,\label{vcres2.17}\\ 
&X_{2t,9t,29\ell,87\ell,1,2}(29n+10t)=Y_{2t,9t,29\ell,87\ell,1,2}(29n+10t)=0, 
\label{vcres2.20} \\
&X_{2t,3t,11\ell,11\ell,1,2}(11n+4t)=Y_{2t,3t,11\ell,11\ell,1,2}(11n+4t)=0,\label{vcres1.5} \\
&X_{6t,t,19\ell,19\ell,1,2}(19n+4t)=Y_{6t,t,19\ell,19\ell,1,2}(19n+4t)=0.\label{vcres1.20}\\
&X_{3t,4t,19\ell,19\ell,1,3}(19n+17t)=Z_{3t,4t,19\ell,19\ell,1,3}(19n+17t)=0.\label{vcres1.21} \\
&X_{t,2t,13\ell,13\ell,1,2}(13n+10t)=Y_{t,2t,13\ell,13\ell,1,2}(13n+10t)=0,\label{vcres1.9}\\
&X_{t,2t,13\ell,13\ell,1,3}(13n+10t)=Z_{t,2t,13\ell,13\ell,1,3}(13n+10t)=0,\label{vcres1.10}\\
&X_{6t,t,13\ell,13\ell,1,3}(13n+t)=Y_{6t,t,13\ell,13\ell,1,3}(13n+t)=0,\label{vcres1.11}\\
&X_{2t,t,7\ell,7\ell,1,3}(7n+6t)=Y_{2t,t,7\ell,7\ell,1,3}(7n+6t)=0,\label{vcres1.2}\\
&X_{8t,3t,17\ell,17\ell,1,8}(17n+16t)=Y_{8t,3t,17\ell,17\ell,1,8}(17n+16t)=0,\label{vcres1.18}\\
&X_{8t,6t,17\ell,17\ell,1,2}(17n+10t)=Y_{8t,6t,17\ell,17\ell,1,2}(17n+10t)=0,\label{vcres1.15}\\
&X_{8t,t,17\ell,17\ell,1,4}(17n+6t)=Y_{8t,t,17\ell,17\ell,1,4}(17n+6t)=0,\label{vcres1.16}\\
&X_{4t,3t,13\ell,13\ell,1,4}(13n+8t)=Y_{4t,3t,13\ell,13\ell,1,4}(13n+8t)=0,\label{vcres1.12}\\
&X_{5t,6t,11\ell,33\ell,2,5}(11n+9t)=Y_{5t,6t,11\ell,33\ell,2,5}(11n+9t)=0.\label{vcres1.11.25}
\end{align}
\end{theorem}
\begin{theorem}\label{vcthm1} Let $\ell \geq 1$ and $t$ satisfying the condition stated in (\ref{C1}). Then 
\begin{align}
&X_{3t,4t,17\ell,51\ell,2,3}(17n+9t)=Z_{3t,4t,17\ell,51\ell,2,3}(17n+9t)=0,\label{vcres2.13}\\
&X_{4t,9t,17\ell,51\ell,3,2}(17n+15t)=Y_{4t,9t,17\ell,51\ell,3,2}(17n+15t)=0,\label{vcres2.14}\\
&X_{4t,3t,17\ell,51\ell,1,6}(17n+11t)=Y_{4t,3t,17\ell,51\ell,1,6}(17n+11t)=0,\label{vcres2.15}\\
&X_{3t,2t,11\ell,33\ell,2,3}(11n+6t)=Z_{3t,2t,11\ell,33\ell,2,3}(11n+6t)=0,\label{vcres2.5}\\ 
&X_{2t,9t,11\ell,33\ell,3,2}(11n+t)=Y_{2t,9t,11\ell,33\ell,3,2}(11n+t)=0,\label{vcres2.6}\\
&X_{2t,3t,11\ell,33\ell,1,6}(11n+10t)=Y_{2t,3t,11\ell,33\ell,1,6}(11n+10t)=0,\label{vcres2.7}\\
&X_{t,2t,11\ell,33\ell,6,1}(11n+4t)=Z_{t,2t,11\ell,33\ell,6,1}(11n+4t)=0,\label{vcres2.8}\\
&X_{t,2t,5\ell,5\ell,3,3}(5n+2t)=Z_{t,2t,5\ell,5\ell,3,3}(5n+2t)=0,\label{vcres1.1} \\
&X_{3t,2t,7\ell,7\ell,3,3}(7n+t)=Z_{3t,2t,7\ell,7\ell,3,3}(7n+t)=0,\label{vcres1.3} \\
&X_{t,2t,7\ell,7\ell,2,3}(7n+4t)=Z_{t,2t,7\ell,7\ell,2,3}(7n+4t)=0,\label{vcres1.4}\\
&X_{4t,t,11\ell,11\ell,1,6}(11n+5t)=Y_{4t,t,11\ell,11\ell,1,6}(11n+5t)=0,\label{vcres1.6}\\
&X_{3t,4t,11\ell,11\ell,2,3}(11n+9t)=Z_{3t,4t,11\ell,11\ell,2,3}(11n+9t)=0,\label{vcres1.7}\\
&X_{2t,3t,11\ell,11\ell,3,6}(11n+t)=Z_{2t,3t,11\ell,11\ell,3,6}(11n+t)=0,\label{vcres1.8}\\
&X_{6t,4t,13\ell,13\ell,3,3}(13n+2t)=Z_{6t,4t,13\ell,13\ell,3,3}(13n+2t)=0,\label{vcres1.13}\\
&X_{4t,t,13\ell,13\ell,3,4}(13n+8t)=Y_{4t,t,13\ell,13\ell,3,4}(13n+8t)=0,\label{vcres1.14}\\
&X_{t,4t,17\ell,17\ell,3,3}(17n+16t)=Z_{t,4t,17\ell,17\ell,3,3}(17n+16t)=0,\label{vcres1.17}\\
&X_{4t,3t,17\ell,17\ell,3,6}(17n+15t)=Y_{4t,3t,17\ell,17\ell,3,6}(17n+15t)=0.\label{vcres1.19}
\end{align}
\end{theorem}
\section{Preliminary}
\begin{definition}
The Ramanujan general theta function \cite[p. 34, Eq. 18.1]{C16} is denoted by $f(a,b)$, and is defined as
   \begin{equation}\label{Rdef}
f(a,b):=\sum_{n=-\infty}^{\infty} a^{n(n+1)/2}b^{n(n-1)/2}, \qquad |ab|<1.
   \end{equation}
\end{definition}
\begin{lemma}{\cite[p. 34, Entry 18 and Entry 19]{C16}} We have
\begin{align}
f(a,b)&=f(b,a), \label{fprop1.0} \\ 
f(1,a)&=2f(a,a^3), \label{fprop1.1} \\
f(-1,a)&=0, \label{fprop1.2} \\
f(a,b)&=(-a,-b,ab;ab)_{\infty}. \label{fprop1.3}  
\end{align}
\end{lemma}

\begin{definition}
    Ramanujan defines the following three particular cases of \eqref{Rdef} \cite[p.36, Entry 22]{C16}:
\begin{align*}
\phi(q)&:=f(q,q),  \\
\varphi(q)&:=f(q,q^3)=\frac{1}{2}f(1,q), \\
f(-q)&:=f(-q,-q^2)=\p(q;q).
\end{align*}
\end{definition}
The following lemma expresses the product of two theta functions as the sum of the product of theta functions:
\begin{lemma}{\cite[p. 46, Entry 30]{C16}} If $ab=cd$, then
\begin{align}
\label{R1}f(c,d)f(a,b)&=af\left({b\over d},{d\over b}abcd \right)f\left({b\over c},{c\over b}abcd \right)+f(ad,bc)f(ac,bd),\\
\label{R2}f^2(a,b)&=af(1,a^2b^2)f(b/a,a^3b)+ f(ab,ab)f(a^2,b^2), \\
\label{R3}f(a,b)&=af(b/a,a^5b^3)+f(a^3b,ab^3).
\end{align}
\end{lemma}
\begin{lemma}\label{f3lemma}\cite[Lemma 2.2]{tang2023vanishing} We have
\begin{nalign}
\label{F1}f(q^k, q^{\mu - k})^3 = &f(q^{3k}, q^{3\mu - 3k}) M(\mu,1) + q^k f(q^{\mu+3k},q^{2\mu-3k})M(\mu,2) \nonumber \\ & + q^{2k}f(q^{\mu-3k},q^{2\mu+3k})M(\mu,2), \\
\label{F2}f(-q^k, -q^{\mu - k})^3 =&f(-q^{3k},-q^{3\mu - 3k}) M(\mu,1)-q^k f(-q^{\mu+3k},-q^{2\mu-3k})M(\mu,2) \nonumber \\ &+ q^{2k}f(-q^{\mu-3k},-q^{2\mu+3k})M(\mu,2)
\end{nalign}
where
\begin{align*}
M(\mu,1)=&f(q^{\mu},q^{\mu})f(q^{3\mu},q^{3\mu}) +  q^{\mu} f(1,q^{2\mu})f(1,q^{6\mu}), \\
M(\mu,2)=&f(1,q^{2\mu})f(q^{2\mu},q^{4\mu}) + f(q^{\mu},q^{\mu})f(q^{\mu},q^{5\mu})  
 \end{align*}
 \end{lemma}
It is worth noting that $M(\mu,1)$ and $M(\mu,2)$ are series in $q^\mu$.
\begin{lemma}{\cite[Lemma 2.2]{liu2024vanishing}}
For any positive integer $n$, we have
\begin{align}\label{nthpower}
f(a, b)^n=\sum_{z=-\infty}^{\infty} C_{n, z}(ab) a^z f\left(a^{n+z} b^z, a^{-z} b^{n-z}\right)
\end{align}
where $C_{n, z}(ab)$ is a formal Laurent series in $ab$.
\end{lemma}
We end this section with the following definition: 
\begin{definition}
Let $k>0, l \geq 0$ be integers and let $A(q)=\sum_{n=0}^{\infty} a(n) q^n$ be a formal power series. Define $E_{k, l}$ by
$$
E_{k, l}(A(q)):=\sum_{n=0}^{\infty} a(k n+l) q^{k n+l} .
$$
\end{definition}

\section{Proofs of \Cref{vcthm1,vcthm2}}
\begin{proof}[Proof of \cref{vcres2.0}]
Upon using \cref{nthpower} with $n=2$, we write
\begin{align}\label{E1}
\sum_{n=-\infty}^{\infty}X_{t,2t,5\ell,15\ell,2,1}(n)q^n=&\p(q^{t},q^{5\ell-t};q^{5\ell})^2\p(q^{2t},q^{15\ell-2t};q^{15\ell}) \nonumber\\
=&\frac{1}{\p(q^{5\ell};q^{5\ell})^2\p(q^{15\ell};q^{15\ell})}\Bigg\{f(-q^t,-q^{5\ell-t})^2f(-q^{2t},-q^{15\ell-2t})\Bigg\} \nonumber\\
=&\frac{1}{\p(q^{5\ell};q^{5\ell})^2\p(q^{15\ell};q^{15\ell})}\Bigg\{\sum_{z=-\infty}^{\infty}(-1)^zC_{2,z}(q^{5\ell})\nonumber \\ &\times q^{tz} f(-q^{2t},-q^{15\ell-2t}) f(q^{5\ell z+2t},q^{-5\ell z+10\ell-2t}) \Bigg\}.
\end{align}
Observe that 
\begin{nalign*}
q^{tz}f(-q^{2t},-q^{15\ell-2t})f(q^{5\ell z+2t},q^{-5\ell z+10\ell-2t})& =\sum_{m,n=-\infty}^{\infty}(-1)^mq^{\frac {15\ell m^{2}}{2}+\frac{15\ell m}{2}-2tm+5\ell n^{2}-5\ell nz
+5\ell n-2nt+tz}.
\end{nalign*}
Employing the linear transformations $m=2r-z+s-1$, $n=3r-z-s$ in above sum, we obtain 
\begin{align*}E_{5,2t}\left(q^{tz}f(-q^{2t},-q^{15\ell-2t})f(q^{5\ell z+2t},q^{-5\ell z+10\ell-2t})\right)=0, \forall ~z.
\end{align*}
Therefore, from (\ref{E1}), we arrive at $X_{t,2t,5\ell,15\ell,2,1}(5n+2t)=0.$
Following a similar strategy, we find $Z_{t,2t,5\ell,15\ell,2,1}(5n+2t)=0$, which establishes the proof of  \cref{vcres2.0}. 
\end{proof}
\begin{proof}[Proof of \cref{vcres1.18}]
We have
\begin{align*}
\sum_{n=-\infty}^{\infty}X_{8t,3t,17\ell,17\ell,1,8}(n)q^n=&\p(q^{8t},q^{17\ell-8t};q^{17\ell})\p(q^{3t},q^{17\ell-3t};q^{17\ell})^8 \\
=&\frac{1}{\p(q^{17\ell};q^{17\ell})^9}\Bigg\{f(-q^{8t},-q^{17\ell-8t})f(-q^{3t},-q^{17\ell-3t})^8\Bigg\}.
\end{align*}
Employing \cref{nthpower} with $n=8$ in the above, we obtain
\begin{nalign*}
\sum_{n=-\infty}^{\infty}X_{8t,3t,17\ell,17\ell,1,8}(n)q^n=
&\frac{1}{\p(q^{17\ell};q^{17\ell})^9}\Bigg\{
\sum_{z=-\infty}^{\infty}(-1)^zC_{8,z}(q^{17\ell})q^{3zt} \\& \times f(-q^{8t},-q^{17\ell-8t})f(q^{17\ell z+24t},q^{-17\ell z+136\ell-24t})\Bigg\}.
\end{nalign*}
To prove \cref{vcres1.18}, it is enough to show that,
\begin{align*}
E_{17,16t}\left(q^{3tz}f(-q^{8t},-q^{17\ell-8t})f(q^{17\ell z+24t},q^{-17\ell z+136\ell-24t})\right)=0, \quad \forall ~z.
\end{align*}
We have
\begin{align*}
~&q^{3tz} f(-q^{8t},-q^{17\ell-8t}) f(q^{17\ell z + 24t},q^{-17\ell z + 136 \ell - 24t}) \\ = ~&  q^{3tz} \sum_{m,n=-\infty}^{\infty} (-1)^m q^{8tm + \frac{17}{2}\ell m^2 - \frac{17}{2} \ell m + 17\ell nz + 24nt + 68\ell n^2 - 68\ell n}.
\end{align*}
In order to extract those terms whose exponent is congruent to $16t \pmod{17}$ in the above, set $8mt+7nt+3tz\equiv 16t \pmod{17}$, which is equivalent to $m+3n\equiv -2-6z \pmod{17}$ and $3m-8n\equiv -6-z \pmod{17}$. Solving these congruences, we obtain $m=8r-3z+3s -2$ and $n=-3r+z+s$. So, $(17n+16t)$-component of $R$ is 
 \begin{nalign*}
&E_{17,16t}(q^{3tz} f(-q^{8t},-q^{17\ell-8t}) f(q^{17\ell z + 24t},q^{-17\ell z + 136 \ell - 24t}))\\=&  \sum_{r,s=-\infty}^{\infty} (-1)^{z+s}q^{17\ell + \frac{17}{2} \ell z  + 16t - 867\ell zr + \frac{289}{2}\ell s^2 - 136tr + 51tz - \frac{289}{2}\ell s + \frac{323}{2}\ell z^2 + 1156\ell r^2} 
     \\ = &(-1)^zq^{17\ell+\frac{17\ell z}{2}+\frac{323\ell z^2}{2}+51tz+16t} f(-1,-q^{289\ell}) f(q^{867\ell z+1156\ell+136t},q^{-867\ell z+1156\ell-136t}) = 0.
 \end{nalign*}
Using the same approach, one can demonstrate $Y_{8t,3t,17\ell,17\ell,1,8}(17n+16t)=0$, which completes the proof of \cref{vcres1.18}. The proofs of \crefrange{vcres2.1}{vcres1.2}, \crefrange{vcres1.15}{vcres1.12} follow the same path. 

\end{proof}

\begin{proof}[Proof of \cref{vcres1.11.25}] From (\ref{R2}) and  \cref{nthpower}, we have
\begin{nalign*}
\sum_{n=-\infty}^{\infty}&X_{5t,6t,11\ell,33\ell,2,5}(n)q^n\\
=&\p(q^{5t},q^{11\ell-5t};q^{11\ell})^2\p(q^{6t},q^{33\ell-6t};q^{33\ell})^5 \\
=&\frac{1}{\p(q^{11\ell};q^{11\ell})^2\p(q^{33\ell};q^{33\ell})^5}\Bigg\{f(-q^{5t},-q^{11\ell-5t})^2f(-q^{6t},-q^{33\ell-6t})^5\Bigg\} \\
=&\frac{1}{\p(q^{11\ell};q^{11\ell})^2\p(q^{33\ell};q^{33\ell})^5}\bigg(\varphi(q^{11\ell})f(q^{10t},q^{22\ell-10t}) -2q^{5t}\psi(q^{22\ell})\\ &f(q^{11\ell-10t},q^{10t+11\ell})\bigg)\times \bigg(\sum_{z=-\infty}^{\infty}(-1)^zC_{5,z}(q^{33\ell})f(-q^{165\ell-33\ell z-30t},-q^{33\ell z+30t})\bigg) \\
=&\frac{1}{\p(q^{11\ell};q^{11\ell})^2\p(q^{33\ell};q^{33\ell})^5}\sum_{z=-\infty}^{\infty}(-1)^zC_{5,z}(q^{33\ell})\bigg(\varphi(q^{11\ell})R_{1,1}-2\psi(q^{22\ell})R_{1,2}\bigg),
\end{nalign*}
where \begin{align}
\label{P1}R_{1,1}:=&q^{6tz}\sum_{m,n=-\infty}^{\infty}(-1)^nq^{11\ell m^{2}+11\ell m-10tm+\frac{165\ell n^{2}}{2}-33\ell nz+\frac{165\ell n}{2}-30nt}, \\
\label{P2}R_{1,2}:=&q^{6tz+5t}\sum_{m,n=-\infty}^{\infty}(-1)^nq^{11\ell m^{2}+10tm+\frac{165\ell n^{2}}{2}-33\ell nz+\frac{165\ell n}{2}-30nt}.
\end{align}
To obtain the $(11n+9t)$-component of $R_{1,1}$, we set $mt+3nt+6zt \equiv 9t \pmod{11}$, which yields $m+3n\equiv -2-6z \pmod{11}$ and $4m-10n\equiv 3-2z \pmod{11}$. Note that $s$ takes the form $s=2k$ for some integer $k$.  Solving  $m+3n= 11r-2-6z$ and $4m-10n=11s-8-2z$ for $m $ and $n$ and then substituting in (\ref{P1}), we obtain 
\begin{nalign*}
E_{11,9t}(R_{1,1})=(-1)^zq^{22\ell +20t+\frac{429}{2}\ell z^2 + 66tz+ \frac{33}{2}\ell z}f(-1,-q^{363\ell})f(q^{-726\ell z+605\ell-110t},q^{726\ell z+605\ell +110t})=0. 
\end{nalign*}
In $R_{1,2}$, settings $10mt+3nt+6tz+5t\equiv 9t \pmod{11}$ gives us
$m-3n=11r-4+6z$ and $4m+10n=11s+6+2z$. From this, it is clear that $s$ takes the form $s=2k$ for integer $k$ and upon solving, we obtain $m=5r-1+3z+3s$, $n=-2r+1-z+s$. Substituting these in \eqref{P2}, we obtain
\begin{nalign*}
 E_{11,9t}(R_{1,2})=(-1)^{1-z}q^{176\ell-35t+\frac{429}{2}\ell z^2+66tz-\frac{693}{2}\ell z}f(-1,-q^{363\ell})f(q^{726\ell z + 110t},q^{-726\ell z + 1210\ell-110t})=0.    
\end{nalign*}
This completes the proof.
\end{proof}
\begin{proof}[Proof of \eqref{vcres2.7}] Using \eqref{F2}, we find that
\begin{nalign*}
\sum_{n=-\infty}^{\infty}&X_{2t,3t,11\ell,33\ell,1,6}(n)q^n \\ 
=&\p(q^{2t},q^{11\ell-2t};q^{11\ell})\p(q^{3t},q^{33\ell-3t};q^{33\ell})^6 \\
=&\frac{1}{\p(q^{11\ell};q^{11\ell})\p(q^{33\ell};q^{33\ell})^6}\Bigg\{f(-q^{2t},-q^{11\ell-2t})f(-q^{3t},-q^{33\ell-3t})^6\Bigg\} \\
=&\frac{1}{\p(q^{11\ell};q^{11\ell})\p(q^{33\ell};q^{33\ell})^6}\Bigg\{f(-q^{2t},-q^{11\ell-2t})\bigg(f(-q^{3t},-q^{33\ell-3t})^3\bigg)^2\Bigg\} \\
=&\frac{1}{\p(q^{11\ell};q^{11\ell})\p(q^{33\ell};q^{33\ell})^6}\Bigg\{f(-q^{2t},-q^{11\ell-2t})  \bigg(f(-q^{9t},-q^{99\ell-9t})M(33\ell,1) \\ &-q^{3t}f(-q^{33\ell+9t},-q^{66\ell-9t})M(33\ell,2) +q^{6t}f(-q^{33\ell-9t},-q^{66\ell+9t}) M(33\ell,2)\bigg)^2\Bigg\}
\end{nalign*}   
which is equivalent to
\begin{nalign*}
\sum_{n=-\infty}^{\infty}& X_{2t,3t,11\ell,33\ell,1,6}(n)q^n \\ =&\frac{1}{\p(q^{11\ell};q^{11\ell})\p(q^{33\ell};q^{33\ell})^6}\Bigg\{\varphi(q^{99\ell})M(33\ell,1)^2R_{2,1} -2\varphi(q^{198\ell})M(33\ell,1)^2R_{2,2} \\
&+2M(33\ell,1)M(33\ell,2)\bigg(f(q^{33\ell},q^{165\ell})\bigg(R_{2,3}+R_{2,4}\bigg) -f(q^{66\ell},q^{132\ell})\bigg(R_{2,5}+R_{2,6}\bigg)\bigg) \\ & +\phi(q^{99\ell})M(33\ell,2)^2\bigg(R_{2,4}-R_{2,3}\bigg)-2q^{33\ell}\varphi(q^{198\ell})M(33\ell,2)^2\bigg(R_{2,6}-R_{2,5}\bigg)\\
&+2q^{33\ell}M(33\ell,2)^2\bigg(f(q^{33\ell},q^{165\ell})R_{2,1}-f(q^{66\ell},q^{132\ell})R_{2,2}\bigg)\Bigg\} \numberthis \label{proofmod1116s1}
\end{nalign*}
where,
\begin{nalign*}
R_{2,1}&:=f(q^{18t},q^{198\ell-18t})f(-q^{2t},-q^{11\ell-2t}) \\& =\sum_{m,n=-\infty}^{\infty}(-1)^nq^{99\ell m^{2}+99\ell m-18tm+\frac {11\ell n^{2}}{2}+\frac{11\ell n}{2}-2tn}, \\
R_{2,2}&:=q^{9t}f(q^{99\ell-18t},q^{99\ell+18t})f(-q^{2t},-q^{11\ell-2t}) \\ & =q^{9t}\sum_{m,n=-\infty}^{\infty}(-1)^nq^{99\ell m^{2}+18tm+\frac {11\ell n^{2}}{2}+\frac{11\ell n}{2}-2tn}, \\
R_{2,3}&:=q^{12t}f(q^{66\ell-18t},q^{132\ell+18t})f(-q^{2t},-q^{11\ell-2t}) \\ & =q^{12t}\sum_{m,n=-\infty}^{\infty}(-1)^nq^{99\ell m^{2}+33\ell m+18tm+\frac{11\ell n^{2}}{2}+\frac {11\ell n}{2
}-2tn}, \\
R_{2,4}&:=q^{6t}f(q^{66\ell+18t},q^{132\ell-18t})f(-q^{2t},-q^{11\ell-2t}) \\ & =q^{6t}\sum_{m,n=-\infty}^{\infty}(-1)^nq^{99\ell m^{2}+33\ell m-18tm+\frac{11\ell n^{2}}{2}+\frac {11\ell n}{2
}-2tn},\\
R_{2,5}&:=q^{3t}f(q^{33\ell+18t},q^{165\ell-18t})f(-q^{2t},-q^{11\ell-2t}) \\ & =q^{3t}\sum_{m,n=-\infty}^{\infty}(-1)^nq^{99\ell m^{2}+66\ell m-18tm+\frac{11\ell n^{2}}{2}+\frac{11\ell n}{2
}-2tn}, \\
R_{2,6}&:=q^{15t}f(q^{33\ell-18t},q^{165\ell+18t})f(-q^{2t},-q^{11\ell-2t}) \\ & =q^{15t}\sum_{m,n=-\infty}^{\infty}(-1)^nq^{99\ell m^{2}+66\ell m+18tm+\frac{11\ell n^{2}}{2}+\frac{11\ell n}{2
}-2tn}.
\end{nalign*}
To arrive at the required result, we extract those terms whose exponents are congruent to $10t \pmod{11}$ in $R_{2, i}$, $1\leq i\leq6$. Substituting  $m=-1+r+s$, $n=4-9r+2s$ in $R_{2,1}$, $m=-1+r+s$, $n=-4+9r-2s$ in $R_{2,2}$, $m=r+$s, $n=9r-2s+1$ in $R_{2,3}$, $m=r+s$, $n=-9r+2s-2$ in $R_{2,4}$,  $m=r+s$, $n=-9r+2s+2$ in $R_{2,5}$  and  $m=r+s$, $n=9r-2s-3$ in $R_{2,6}$, we find that

\begin{align*}
E_{11,10t}(R_{2,1})&=\sum_{r,s=-\infty}^{\infty}(-1)^rq^{\frac {1089\ell r^{2}}{2}+121\ell s^{2}-\frac {1089\ell r}{2}+110\,l-22ts+10t} \\
&=q^{110\ell+10t}f(-1,-q^{1089\ell})f(q^{121\ell-22t},q^{121\ell+22t})=0,\\
E_{11,10t}(R_{2,2})=&\sum_{r,s=-\infty}^{\infty}(-1)^rq^{\frac{1089\ell r^{2}}{2}+121\ell {s}^{2}-\frac {1089\ell r}{2}-121\ell s+165\ell+22ts-t} \\
&=q^{165\ell-t}f(-1,-q^{1089\ell})f(q^{22t},q^{242\ell-22t})=0, \\
E_{11,10t}(R_{2,3})=&\sum_{r,s=-\infty}^{\infty}(-1)^{1+r}q^{\frac{1089\ell r^{2}}{2}+121\ell s^{2}+\frac {363\ell r}{2}+22ts+11\ell+10t} \\
&=-q^{11\ell+10t}f(-q^{726\ell},-q^{363\ell})f(q^{121\ell-22t},q^{121\ell+22t}),\\
E_{11,10t}(R_{2,4})=&\sum_{r,s=-\infty}^{\infty}(-1)^{r}q^{\frac{1089\ell r^{2}}{2}+121\ell s^{2}+\frac {363\ell r}{2}-22ts+11\ell+10t} \\
&=q^{11\ell+10t}f(-q^{726\ell},-q^{363\ell})f(q^{121\ell-22t},q^{121\ell+22t}), \\
E_{11,10t}(R_{2,5})=&\sum_{r,s=-\infty}^{\infty}(-1)^{r}q^{\frac{1089\ell r^{2}}{2}+121\ell s^{2}-\frac {363\ell r}{2}+121\ell s-22ts+33\ell-t} \\
&=q^{33\ell-t}f(-q^{726\ell},-q^{363\ell})f(q^{242\ell-22t},q^{22t}), \\
E_{11,10t}(R_{2,6})=&\sum_{r,s=-\infty}^{\infty}(-1)^{1+r}q^{\frac {1089\ell r^{2}}{2}+121\ell s^{2}-\frac{363\ell r}{2}+121\ell s+22ts+33\ell+21t} \\
&=-q^{33\ell+21t}f(-q^{726\ell},-q^{363\ell})f(q^{242\ell+22t},q^{-22t}),\\&=-q^{33\ell-t}f(-q^{726\ell},-q^{363\ell})f(q^{242\ell-22t},q^{22t}).
\end{align*}
Using these in \cref{proofmod1116s1}, we deduce that $X_{2t,3t,11\ell,33\ell,1,6}(11n+10t)=0$. In a similar way, one can prove $Y_{2t,3t,11\ell,33\ell,1,6}(11n+10t)=0$. This complete the proof of \cref{vcres2.7}. The proofs of \eqref{vcres2.15}, \eqref{vcres1.6} and \eqref{vcres2.8} follow in the same path. We omit the details here. 
\end{proof}
\begin{proof}[Proof of \eqref{vcres1.4}] Applying (\ref{R2}) and (\ref{F2}),  we obtain
\begin{nalign*}
\sum_{n=-\infty}^{\infty}X_{t,2t,7\ell,7\ell,2,3}(n)q^n=&\p(q^{t},q^{7\ell-t};q^{7\ell})^2\p(q^{2t},q^{7\ell-2t};q^{7\ell})^3 \\
=&\frac{1}{\p(q^{7\ell};q^{7\ell})^5}\Bigg\{f(-q^{t},-q^{7\ell-t})^2f(-q^{2t},-q^{7\ell-2t})^3\Bigg\} \\
=&\frac{1}{\p(q^{7\ell};q^{7\ell})^5}\Bigg\{\bigg(f(q^{2t},q^{14\ell-2t}) \varphi(q^{7\ell})-2q^{t}f(q^{7\ell-2t},q^{2t+7\ell})\varphi(q^{14\ell})\bigg)\\
&\times \bigg(f(-q^{6t},-q^{21\ell-6t})M(7\ell,1)-q^{2t}f(-q^{7\ell+6t},-q^{14\ell-6t})M(7\ell,2) \\ &+q^{4t}f(-q^{7\ell-6t},-q^{14\ell+6t})M(7\ell,2)
\bigg)\Bigg\} \\
=&\frac{1}{\p(q^{7\ell};q^{7\ell})^5}\Bigg\{\varphi(q^{7\ell})M(7\ell,1)R_{3,1}-2\varphi(q^{14\ell})M(7\ell,1)R_{3,2} \\
&+\varphi(q^{7\ell})M(7\ell,2)\bigg(R_{3,3}-R_{3,4}\bigg)+2\varphi(q^{14\ell})M(7\ell,2)\bigg(R_{3,5}-R_{3,6}\bigg)\Bigg\}, \numberthis \label{proof1.1s1}
\end{nalign*}
where
\begin{nalign*}
R_{3,1}&:f(q^{14\ell-2t},q^{2t}) f(-q^{21\ell-6t},-q^{6t}) =\sum_{m,n=-\infty}^{\infty}(-1)^nq^{7\ell m^{2}+7\ell m-2tm+{\frac{21\ell n^{2}}{2}}+{\frac {21\ell n}{2}}-
6tn}\\
R_{3,2}&:=q^t f(q^{7\ell+2t},q^{7\ell-2t}) f(-q^{21\ell-6t},-q^{6t})=q^t\sum_{m,n=-\infty}^{\infty}(-1)^nq^{7\ell m^{2}+2tm+\frac {21\ell n^{2}}{2}+\frac {21\ell n}{2}-6tn}, \\
R_{3,3}&:=q^{4t} f(q^{14\ell-2t},q^{2t}) f(-q^{14\ell+6t},-q^{7\ell-6t}) = q^{4t}\sum_{m,n=-\infty}^{\infty}(-1)^nq^{7\ell m^{2}+7\ell m-2tm+\frac{21\ell n^{2}}{2}+\frac{7\ell n}{2}+6tn}\\ 
R_{3,4}&:=q^{2t} f(q^{14\ell-2t},q^{2t}) f(-q^{14\ell-6t},-q^{7\ell+6t}) = q^{2t}\sum_{m,n=-\infty}^{\infty}(-1)^nq^{7\ell m^{2}+7\ell m-2tm+\frac{21\ell n^{2}}{2}+\frac{7\ell n}{2}-6tn}, \\
R_{3,5}&:=q^{3t} f(q^{7\ell+2t},q^{7\ell-2t}) f(-q^{14\ell-6t},-q^{7\ell+6t}) =q^{3t}\sum_{m,n=-\infty}^{\infty}(-1)^nq^{7\ell m^{2}+2tm+\frac {21\,l{n}^{2}}{2}+\frac {7\ell n}{2}-6tn}\\ 
R_{3,6}&:=q^{5t} f(q^{7\ell+2t},q^{7\ell-2t}) f(-q^{14\ell-6t},-q^{7\ell+6t}) = q^{5t}\sum_{m,n=-\infty}^{\infty}(-1)^nq^{7\ell m^{2}+2tm+\frac {21\ell{n}^{2}}{2}+\frac {7\ell n}{2}-6tn}. 
\end{nalign*}
To arrive at the required result, we extract those terms whose exponents are congruent to $4t \pmod{7}$. Employ the transformations: $m=r+2+3s$, $n=2r+1-s$ in $R_{3,1}$, $m=r+2+3s$, $n=-2r-1+s$ in $R_{3,2}$, $m=r+3s$, $n=-2r+s$ in $R_{3,3}$, $m=r+3s-1$, $n=2r-s$ in $R_{3,4}$, $m=r+3s$, $n=-2r+1+s$ in $R_{3,5}$, and finally $m=r+2+3s$, $n=2r+1-s$ in $R_{3,6}$. We obtain $(7n+4t)$-component
\begin{align*}
E_{7,4t}(R_{3,1})=&E_{7,4t}(R_{3,2})=0, \\
E_{7,4t}(R_{3,3})=&E_{7,4t}(R_{3,4})=q^{4t}f(q^{49\ell+14t},q^{49\ell-14t})f(-q^{49\ell},-q^{98\ell}),\\
E_{7,4t}(R_{3,5})=&E_{7,4t}(R_{3,6})=- q^{14\ell-3t}f(q^{98\ell-14t},q^{14t})f(-q^{49\ell},-q^{98\ell}).
\end{align*}
Operating $E_{7,4t}$ to \eqref{proof1.1s1} on both sides and employing the above, we arrive at $X_{t,2t,7\ell,7\ell,2,3}(7n+4t)=0$. Similarly, we obtain $Z_{t,2t,7\ell,7\ell,2,3}(7n+4t)=0$, completing the proof of \eqref{vcres1.4}. The proofs of 
\eqref{vcres2.13}, \eqref{vcres2.14}, \eqref{vcres2.5}, \eqref{vcres2.6} and \eqref{vcres1.7} follow in the same path.
\end{proof}
\begin{proof}[Proof of \eqref{vcres1.14}] From (\ref{R2}) and (\ref{F2}), we have
\begin{nalign*}
 &\sum_{n=-\infty}^{\infty}X_{4t,t,13\ell,13\ell,3,4}(n)q^n\\=& \p(q^{4t},q^{13\ell-4t};q^{13\ell})^3\p(q^{t},q^{13\ell-t};q^{33\ell})^4 
= \frac{1}{\p(q^{13\ell};q^{13\ell})^7}  \Bigg \{f^3(-q^{4t},-q^{13\ell-4t}) f^4(-q^t, -q^{13\ell-t}) \Bigg \}\\
= & \frac{1}{\p(q^{13\ell},q^{13\ell})^7}  \Bigg \{ \bigg \{ f(-q^{12t},-q^{39\ell-12t}M(13\ell,1)  - q^{4t} f(-q^{13\ell+12t},-q^{26\ell-12t})M(13\ell,2) \\ &+ q^{8t} f(-q^{13\ell+12t},-q^{26\ell+12t})M(13\ell,2) \bigg \}  \bigg (f^2(-q^t,-q^{13\ell-t}) \bigg )^2 \Bigg \}\\
= & \frac{1}{\p(q^{13\ell};q^{13\ell})^7} \Bigg \{ R_{6,1} M(13\ell,1) \bigg ( f^2(1,q^{26\ell}) f(1,q^{52\ell}) q^{13\ell} + \varphi(q^{26\ell}) \varphi^2(q^{13\ell}) \bigg )   + R_{6,2} M(13\ell,1)\\ &\times \bigg ( f(1,q^{52\ell}) \varphi^2(q^{13\ell}) + f^2(1,q^{26\ell}) \varphi(q^{26\ell}) \bigg ) -  2 M(13\ell,1) \varphi(q^{13\ell}) f(1,q^{26\ell}) f(q^{13\ell},q^{39\ell}) \\ & \times \bigg ( R_{6,3} + R_{6,4} \bigg )  M(13\ell,2) \Bigg \{ - \bigg \{ \varphi^2(q^{13\ell}) \bigg ( \varphi(q^{26\ell}) \bigg ( R_{6,5} - R_{6,6} \bigg ) + f(1,q^{53\ell}) \\ & \times \bigg ( R_{6,7} - R_{6,8} \bigg ) \bigg ) \bigg \}  + \bigg \{ 2 \varphi(q^{13\ell}) f(1,q^{26\ell}) f(q^{13\ell},q^{39\ell}) \bigg ( R_{6,9} - R_{6,10} + R_{6,11} - R_{6,12} \bigg ) \bigg \} \\ & - \bigg \{ f^2(1,q^{26\ell}) \bigg ( \varphi (q^{26\ell}) \bigg ( R_{6,6} - R_{6,10} \bigg ) + f(1,q^{52\ell}) q^{13\ell} \bigg ( R_{6,5} - R_{6,9} \bigg ) \bigg ) \bigg \} \Bigg \}, \numberthis \label{proof34}
\end{nalign*}
where
\begin{align*}
&R_{6,1} := 
\sum_{m,n=-\infty}^{\infty} q^{12tm+\frac{39}{2}\ell m^2 -\frac{39}{2}\ell m+4tn+26\ell n^2-26\ell n},\\&
R_{6,2}:= 
q^{2t}  \sum_{m,n=-\infty}^{\infty} (-1)^m q^{12tm+\frac{39}{2}\ell m^2 -\frac{39}{2}\ell m-4tn+26\ell n^2},\\
&R_{6,3}:= 
q^{t}  \sum_{m,n=-\infty}^{\infty} (-1)^m q^{12tm+\frac{39}{2}\ell m^2 -\frac{39}{2}\ell m+4tn+26\ell n^2-13\ell n}, \\
&R_{6,4}:= 
q^{3t}  \sum_{m,n=-\infty}^{\infty} (-1)^m q^{12tm+\frac{39}{2}\ell m^2 -\frac{39}{2}\ell m-4tn+26\ell n^2-13\ell n},\\
&R_{6,5}:= 
q^{4t}  \sum_{m,n=-\infty}^{\infty} (-1)^m q^{12tm+\frac{39}{2}\ell m^2 -\frac{39}{2}\ell m+4tn+26\ell n^2-26\ell n},\\&
R_{6,6}:= 
q^{8t}  \sum_{m,n=-\infty}^{\infty} (-1)^m q^{-12tm+\frac{39}{2}\ell m^2 -\frac{13}{2}\ell m+4tn+26\ell n^2-26\ell n},\\
&R_{6,7}:= 
q^{6t}  \sum_{m,n=-\infty}^{\infty} (-1)^m q^{12tm+\frac{39}{2}\ell m^2 -\frac{13}{2}\ell m-4tn+26\ell n^2},\\&
R_{6,8}:= 
q^{10t}  \sum_{m,n=-\infty}^{\infty} (-1)^m q^{-12tm+\frac{39}{2}\ell m^2 -\frac{13}{2}\ell m-4tn+26\ell n^2},\\
&R_{6,9}:= 
q^{5t}  \sum_{m,n=-\infty}^{\infty} (-1)^m q^{12tm+\frac{39}{2}\ell m^2 -\frac{13}{2}\ell m+4tn+26\ell n^2-13\ell n},\\&
 R_{6,10} := 
 q^{9t}  \sum_{m,n=-\infty}^{\infty} (-1)^m q^{-12tm+\frac{39}{2}\ell m^2 -\frac{13}{2}\ell m+4tn+26\ell n^2-13\ell n},\\
&R_{6,11} := 
q^{7t}  \sum_{m,n=-\infty}^{\infty} (-1)^m q^{12tm+\frac{39}{2}\ell m^2 -\frac{13}{2}\ell m-4tn+26\ell n^2-13\ell n}, \\&R_{6,12} :=  
q^{11t}  \sum_{m,n=-\infty}^{\infty} (-1)^m q^{-12tm+\frac{39}{2}\ell m^2 -\frac{13}{2}\ell m-4tn+26\ell n^2-13\ell n}.
\end{align*}
To arrive at the required result, we extract those terms in which the exponents are congruent to $8t\pmod{13}$ on both sides of (\ref{proof34}), which can be done using the transformations listed in the following table.
\begin{table}[H]
    \centering
    \begin{tabular}{c|c|c|c}
    \hline
       Term  &  Transformations & Term  &  Transformations \\
       \hline
       $R_{6,1}$  &  $m = 1+4r+s$, $n=-1+r-3s$ &
        $R_{6,2}$  &  $m = -2+4r+s$, $n=-1-r+3s$ \\
       \hline
         $R_{6,3}$  &  $m = 2+4r+s$, $n=-1+r-3s$ &
          $R_{6,4}$  & $m = -1+4r+s$, $n=-1-r+3s$ \\
       \hline
        $R_{6,5}$  &  $m = 4r+s$, $n=1+r-3s$ &
         $R_{6,6}$  &  $m = -4r-s$, $n=r-3s$ \\
       \hline
         $R_{6,7}$  &  $m = -2+4r+s$, $n=-r+3s$ &
        $R_{6,8}$  &  $m = 2+4r-s$, $n=1-r+3s$ \\
       \hline
         $R_{6,9}$  &  $m = 1+4r+s$, $n=1+r-3s$ &
       $R_{6,10}$  & $m = -1-4r-s$, $n=r-3s$ \\
       \hline
       $R_{6,11}$  & $m = -1+4r+s$, $n=-r+3s$ &
       $R_{6,12}$  & $m = 1-4r-s$, $n=1-r+3s$ \\
       \hline
    \end{tabular}
\end{table}
With the help of above table, we obtain the $(13n+8t)$-components as follows.
\begin{align*}
&E_{13,8t}(R_{6,1})=E_{13,8t}(R_{6,2})= E_{13,8t}(R_{6,3}) =  E_{13,8t}(R_{6,4}) =0,\\
&E_{13,8t}(R_{6,5}) = E_{13,8t}(R_{6,6})= q^{8t} f(-q^{169\ell},-q^{338\ell}) f(q^{338\ell+52t},q^{338\ell-52t}),\\
&E_{13,8t}(R_{6,7}) =E_{13,8t}(R_{6,8})=  q^{91\ell-18t} f(-q^{169\ell},-q^{338\ell}) f(q^{52t},q^{676\ell-52t}),\\
&E_{13,8t}(R_{6,9}) = E_{13,8t}(R_{6,10})=-q^{26\ell+21t} f(-q^{169\ell},-q^{338\ell}) f(q^{507\ell+52t},q^{169\ell-52t}),\\
&E_{13,8t}(R_{6,11}) = E_{13,8t}(R_{6,12})=-q^{26\ell-5t} f(-q^{169\ell},-q^{338\ell}) f(q^{169\ell+52t},q^{507\ell-52t}).
\end{align*}
Employing the above in equation \cref{proof34}, we arrive at the required result. With the similar strategy, one can prove $Y_{4t,t,13\ell,13\ell,3,4}(13n+8t)=0$ completing the proof of \cref{vcres1.14}.
\end{proof}
\begin{proof}[Proof of \cref{vcres1.17}] Using (\ref{F2}), we find that
\begin{nalign}\label{R5}
&\sum_{n=-\infty}^{\infty}X_{t,4t,17\ell,17\ell,3,3}(n) q^n\nonumber\\ &=  ( q^{t}, q^{17\ell-t};q^{17\ell})^3_\infty ( q^{4t}, q^{17\ell-4t};q^{17\ell})^3_\infty=  \frac{1}{(q^{17\ell};q^{17\ell})^6_\infty }\Bigg\{f^3(- q^{t},- q^{17\ell-t}) f^3(- q^{4t}, -q^{17\ell - 4t}) \Bigg \}\nonumber\\ 
  = &  \frac{1}{(q^{17\ell};q^{17\ell})^6_\infty } \Bigg \{ \bigg( f(-q^{3t}, -q^{51\ell -3t}) M(17\ell,1) - q^t f(-q^{17\ell + 3t}, -q^{34\ell - 3t})M(17\ell,2)  + q^{2t}\nonumber\\  &\times f(-q^{17\ell - 3t}, -q^{34\ell + 3t}) M(17\ell,2) \bigg) \bigg( f(-q^{12t}, -q^{51\ell - 12t})M(17\ell,1)  - q^{4t}f(-q^{17\ell + 12t}, -q^{34\ell - 12t} )\nonumber\\  &\times M(17\ell,2) + q^{8t}f(-q^{17\ell - 12t}, -q^{34\ell + 12t})M(17\ell,2) \bigg) \Bigg \}\nonumber \\
  = & \frac{1}{(q^{17\ell};q^{17\ell})^6_\infty } \Bigg \{ R_{5,1} M(17\ell,2)^2 +M(17\ell,2)^2 \bigg(R_{5,5}- R_{5,6} - R_{5,8}+ R_{5,9}  \bigg)+M(17\ell,1) M(17\ell,2)\nonumber\\&\times\bigg( - R_{5,2}+ R_{5,3}- R_{5,4} + R_{5,7}\bigg) \Bigg \},
  \end{nalign}
  where
\begin{nalign*}
R_{5,1} &:
= \sum_{m,n=-\infty}^{\infty} (-1)^{m+n} q^{3tm+\frac{51}{2}\ell m^2-\frac{51}{2}\ell m+12tn+\frac{51}{2}\ell n^2-\frac{51}{2}\ell n},\\
R_{5,2} & :
=  q^{4t} \sum_{m,n=-\infty}^{\infty} (-1)^{m+n} q^{3tm+\frac{51}{2}\ell m^2-\frac{51}{2}\ell m+12tn+\frac{51}{2}\ell n^2-\frac{17}{2}\ell n},\\
R_{5,3} & : 
= q^{8t} \sum_{m,n=-\infty}^{\infty} (-1)^{m+n} q^{3tm+\frac{51}{2}\ell m^2-\frac{51}{2}\ell m-12tn+\frac{51}{2}\ell n^2-\frac{17}{2}\ell n},\\
R_{5,4} & : 
= q^{t} \sum_{m,n=-\infty}^{\infty} (-1)^{m+n} q^{3tm+\frac{51}{2}\ell m^2-\frac{17}{2}\ell m+12tn+\frac{51}{2}\ell n^2-\frac{51}{2}\ell n},\\
R_{5,5} & : 
=q^{5t} \sum_{m,n=-\infty}^{\infty} (-1)^{m+n} q^{3tm+\frac{51}{2}\ell m^2-\frac{17}{2}\ell m+12tn+\frac{51}{2}\ell n^2-\frac{17}{2}\ell n},\\
R_{5,6} & : 
=q^{9t} \sum_{m,n=-\infty}^{\infty} (-1)^{m+n} q^{3tm+\frac{51}{2}lm^2-\frac{17}{2}\ell m-12tn+\frac{51}{2}\ell n^2-\frac{17}{2}\ell n},\\
R_{5,7} & :
= q^{2t} \sum_{m,n=-\infty}^{\infty} (-1)^{m+n} q^{-3tm+\frac{51}{2}\ell m^2-\frac{17}{2}\ell m+12tn+\frac{51}{2}\ell n^2-\frac{51}{2}\ell n},\\
R_{5,8} & :
=q^{6t} \sum_{m,n=-\infty}^{\infty} (-1)^{m+n} q^{-3tm+\frac{51}{2}\ell m^2-\frac{17}{2}\ell m+12tn+\frac{51}{2}\ell n^2-\frac{17}{2}\ell n},\\
R_{5,9} & : 
= q^{10t} \sum_{m,n=-\infty}^{\infty} (-1)^{m+n} q^{-3tm+\frac{51}{2}\ell m^2-\frac{17}{2}\ell m-12tn+\frac{51}{2}\ell n^2-\frac{17}{2}\ell n}.
\end{nalign*}
To arrive at the required result, we extract those terms which are congruent to $16t \pmod{17}$ on both sides of (\ref{R5}). Take $m=-4s+r+2$,$n=-2+4r+s$ in $R_{5,1}$,  $m=r-4s$,$n=1+4r+s$ in $R_{5,2}$, $m=1+r-4s$,$ n=1-4r-s$ in $R_{5,3}$, $m=1+r-4s$, $n=1+4r+s$ in $R_{5,4}$, $m=-2+r-4s$,$n=4r+s$ in $R_{5,5}$, $m=-4r-s$,$n=-2-r+4s$ in $R_{5,6}$, $m=1-r+4s$,$n=4r+s$ in $R_{5,7}$, $m=-1+4r+s$,$n=2+r-4s$ in $R_{5,8}$ and   $m=-2-r+4s$,$n=-4r-s$ in $R_{5,9}$, we find that
\begin{align*}
E_{17,16t}(R_{5,1})=&E_{17,16t}(R_{5,5})=E_{17,16t}(R_{5,9})=0,\\
-E_{17,16t}(R_{5,2})=&E_{17,16t}(R_{5,4})=q^{17\ell+6t} f(-q^{578\ell},-q^{289\ell})f(-q^{578\ell
+51t},-q^{289\ell-51t}),\\
E_{17,16t}(R_{5,3})=&-E_{17,16t}(R_{5,7})=q^{17\ell-t} f(-q^{578\ell},-q^{289\ell})f(-q^{289\ell+51t},-q^{578\ell-51t}),\\
E_{17,16t}(R_{5,6})=&-E_{17,16t}(R_{5,8})=q^{119\ell+33t} f(-q^{578\ell},-q^{289\ell})f(-q^{
-51t},-q^{867\ell+51t}).
\end{align*}
Using these, we deduce that $X_{4t,3t,13\ell,13\ell,1,4}(13n+8t)=0$. With the similar strategy, one can prove $Z_{t,4t,17\ell,17\ell,3,3}(17n+16t)=0$ completing the proof of \cref{vcres1.17}. The proofs of \eqref{vcres1.1}, \eqref{vcres1.3} and \eqref{vcres1.13}  follow the same path. 
\end{proof}
\begin{proof}[Proof of \eqref{vcres1.19}] We have
\begin{nalign*}
\sum_{n=-\infty}^{\infty}&X_{4t,3t,17\ell,17\ell,3,6}(n)q^n \\&=\p(q^{4t},q^{17\ell-4t};q^{17\ell})^3 \p( q^{3t},q^{17\ell-3t};q^{17\ell})^6\\
 & = \frac{1}{(q^{17\ell};q^{17\ell})^9_\infty }\Bigg\{f^3(-q^{4t},-q^{17\ell-4t}) f^6(-q^{3t}, - q^{17\ell - 3t}) \Bigg \}\\
& = \frac{1}{\p (q^{17\ell};q^{17\ell})^9} \Bigg \{ M(17\ell,1)^3 \bigg (  \varphi(q^{51\ell})R_{4,1}- f(1,q^{102\ell})R_{4,2} \bigg )  + M(17\ell,1)^2 M(17\ell,2) \\ & \times \bigg \{ \varphi(q^{51\ell}) \bigg ( -R_{4,3} + R_{4,5} \bigg )  + f(1,q^{102\ell}) \bigg ( R_{4,4} - R_{4,6} \bigg ) - 2 f(q^{34\ell},q^{68\ell}) \bigg (  R_{4,7} + R_{4,14} \bigg ) \\ & +2 f(q^{17\ell},q^{85\ell})  \bigg (  R_{4,8} +  R_{4,13} \bigg ) \bigg \} + M(17\ell,2)^3 \bigg \{ \varphi(q^{51\ell}) \bigg ( -R_{4,15} + R_{4,17} - R_{4,10} \\ & + R_{4,12} \bigg )   + f(1,q^{102\ell}) \bigg ( R_{4,16} - R_{4,18} + R_{4,9} - R_{4,11} \bigg )   +2 f(q^{34\ell},q^{68\ell}) \bigg ( R_{4,4} - R_{4,6} \bigg ) \\ &  + 2 f(q^{17\ell},q^{85\ell}) \bigg ( - R_{4,3} +  R_{4,5} \bigg ) \bigg \}  + M(17\ell,1) M(17\ell,2)^2 \bigg \{ \varphi(q^{51\ell}) \bigg ( R_{4,13} + R_{4,8} \bigg ) \\ & + f(1,q^{102\ell})  \bigg ( -R_{4,14} - R_{4,7} \bigg )  + 2 f(q^{34\ell},q^{68\ell}) \bigg (  R_{4,9} -  R_{4,11} +  R_{4,16} -  R_{4,18} -  R_{4,2} \bigg )  \\&+ 2  f(q^{17\ell},q^{85\ell})\bigg  ( - R_{4,10} +  R_{4,12} -  R_{4,15} +  R_{4,17} + R_{4,1} \bigg )\Bigg\}, \numberthis \label{proof36}
\end{nalign*}
where
\begin{nalign*}
&R_{4,1}:= f(-q^{12t},-q^{51\ell-12t}) f(q^{18t},q^{102\ell-18t}),\\
& R_{4,2} :=q^{9t} f(-q^{12t},-q^{51\ell-12t}) f(q^{51\ell-18t},q^{51\ell+18t}),\\ 
&R_{4,3}:=q^{4t} f(-q^{17\ell+12t},-q^{34\ell-12t})f(q^{18t},q^{102\ell-18t}),\\
&R_{4,4}:=q^{13t} f(-q^{17\ell+12t},-q^{34\ell-12t}) f(q^{51\ell-18t},q^{51\ell+18t}), \\
&R_{4,5}:=q^{8t} f(-q^{17\ell-12t},-q^{34\ell+12t}) f(q^{18t},q^{102\ell-18t}), \\
&R_{4,6}:=q^{17t} f(-q^{17\ell-12t},-q^{34\ell+12t}) f(q^{51\ell-18t},q^{51\ell+18t}),\\ 
&R_{4,7}:=q^{3t} f(-q^{12t},-q^{51\ell-12t})f(q^{17\ell+18t},q^{85\ell-18t}), \\
&R_{4,8}:=q^{12t} f(-q^{12t},-q^{51\ell-12t}) f(q^{34\ell-18t},q^{68\ell+18t}), \\ 
&R_{4,9}:=q^{7t} f(-q^{17\ell+12t},-q^{34\ell-12t}) f(q^{17\ell+18t},q^{85\ell-18t}),\\
&R_{4,10}:=q^{16t} f(-q^{17\ell+12t},-q^{34\ell-12t}) f(q^{34\ell-18t},q^{68\ell+18t}), \\ 
&R_{4,11}:=q^{11t} f(-q^{17\ell-12t},-q^{34\ell+12t}) f(q^{17\ell+18t},q^{85\ell-18t}), \\
&R_{4,12}:=q^{20t} f(-q^{17\ell-12t},-q^{34\ell+12t}) f(q^{34\ell-18t},q^{68\ell+18t}), \\ 
&R_{4,13}: =q^{6t} f(-q^{12t},-q^{51\ell-12t}) f(q^{34\ell+18t},q^{68\ell-18t}), \\
&R_{4,14}:=q^{15t}f(-q^{12t},-q^{51\ell-12t}) f(q^{17\ell-18t},q^{85\ell+18t}), \\
&R_{4,15}:=q^{10t} f(-q^{17\ell+12t},-q^{34\ell-12t})f(q^{34\ell+18t},q^{68\ell-18t}),\\
&R_{4,16}:=q^{19t}f(-q^{17\ell+12t},-q^{34\ell-12t}) f(q^{17\ell-18t},q^{85\ell+18t}), \\ 
&R_{4,17}:=q^{14t} f(-q^{17\ell-12t},-q^{34\ell+12t})f(q^{34\ell+18t},q^{68\ell-18t}), \\
&R_{4,18}:=q^{23t} f(-q^{17\ell-12t},-q^{34\ell+12t}) f(q^{17\ell-18t},q^{85\ell+18t}).
\end{nalign*}
To arrive at the required result, we pick those terms whose exponents are congruent $15t \pmod{17}$ on both sides of (\ref{proof36}). The linear transformations to be used to extract the $(17n+15t)$-component for each of the $R_{4,i}$'s are summarized the following table:
    \begin{table}[H]
        \centering
        \begin{tabular}{c|c|c|c}
        \hline
            Term & Transformations &Term & Transformations \\
            \hline
             $R_{4,1}$ &   $m=4r+3s$, $n = -2+3r-2s$ 
            &$R_{4,2}$ &   $m=-1+4r+3s$, $n = -1-3r+2s$ \\
              \hline
               $R_{4,3}$ &  $m=1+4r+3s$, $n = -1-3r-2s$ &
               $R_{4,4}$ &   $m=4r+3s$, $n = -2-3r+2s$  \\
              \hline
               $R_{4,5}$ &   $m=-2-4r-3s$, $n = 3r-2s$  &
              $R_{4,6}$ &   $m=-4r-3s$, $n =2- 3r+2s$  \\
              \hline
              $R_{4,7}$ &   $m=1+4r+3s$, $n =3r-2s$ &
              $R_{4,8}$ &  $m=3+4r+3s$, $n =-1-3r+2s$\\
              \hline
              $R_{4,9}$ &   $m=2+4r+3s$, $n =1+3r-2s$ &
               $R_{4,10}$ &  $m=4r+3s$, $n =1-3r+2s$ \\
              \hline
               $R_{4,11}$ &   $m=1-4r-3s$, $n =-1+3r-2s$  &
              $R_{4,12}$ &   $m=-1-4r-3s$, $n =-3r+2s$\\
              \hline
               $R_{4,13}$ &   $m=2+4r+3s$, $n =2+3r-2s$ &
               $R_{4,14}$ &   $m=4r+3s$, $n =-3r+2s$\\
              \hline
              $R_{4,15}$ &  $m=-1+4r+3s$, $n =3r-2s$ &
               $R_{4,16}$ &   $m=1+4r+3s$, $n =-1-3r+2s$\\
              \hline
               $R_{4,17}$ & $m=-4r-3s$, $n =1-3r-2s$ &
              $R_{4,18}$ &  $m=2-4r-3s$, $n =1-3r+2s$ \\
              \hline
        \end{tabular}
    \end{table} 
Using the above table, we obtain $(17n+15t)$-component for each of the $R_{4,i}$, $1 \leq i \leq 18$ as below.
\begin{nalign*}
&E_{17,15t}(R_{4,1})=E_{17,15t}(R_{4,2})=E_{17,15t}(R_{4,3})=
E_{17,15t}(R_{4,4})=E_{17,15t}(R_{4,5})=E_{17,15t}(R_{4,6})=0,\\
-&E_{17,15t}(R_{4,7})=E_{17,15t}(R_{4,14})=q^{15t}f(-q^{578\ell},-q^{289\ell}) f(q^{867\ell+102t},q^{867\ell-102t}) ,\\
-&E_{17,15t}(R_{4,8}) =E_{17,15t}(R_{4,13})=q^{221\ell+66t}f(-q^{578\ell},-q^{289\ell}) f(q^{-102t},q^{1734\ell+102t}) ,\\
&E_{17,15t}(R_{4,9}) =  -E_{17,15t}(R_{4,16})=q^{102\ell+49t}f(-q^{578\ell},-q^{289\ell}) f(q^{1445\ell+102t},q^{289\ell-102t}) ,\\
&E_{17,15t}(R_{4,10}) = -E_{17,15t}(R_{4,15})= q^{34\ell-2t}f(-q^{578\ell},-q^{289\ell}) f(q^{578\ell+102t},q^{1156\ell-102t}) ,\\
-&E_{17,15t}(R_{4,11}) =E_{17,15t}(R_{4,18})=q^{102\ell-19t}f(-q^{578\ell},-q^{289\ell}) f(q^{289\ell+102t},q^{1445\ell-102t}) ,\\
-&E_{17,15t}(R_{4,12}) =E_{17,15t}(R_{4,17})=q^{34\ell+32t}f(-q^{578\ell},-q^{289\ell}) f(q^{1156\ell+102t},q^{578\ell-102t}).
\end{nalign*}  
Using these in \cref{proof36}, we deduce that $X_{4t,3t,17\ell,17\ell,3,6}(5n+2t)=0$. With the similar strategy, one can prove $Y_{4t,3t,17\ell,17\ell,3,6}(17n+15t)=0$ completing the proof of \cref{vcres1.19}. The proof of \eqref{vcres1.8} follows in the same path.
\end{proof}

\section{Conclusions}
 In this paper, we proved new set of families of vanishing coefficient results for (\ref{D1})-(\ref{D4}) by using Ramanujan's theta function and Jacobi Triple Product Identity. It would be interesting to find what would be the results for other values of $s$,$k$,$u$ and $v$, which are not considered in this paper.
\section*{Declarations }
\textbf{Conflict of interest}: The authors declares no conflict of interest.
\bibliographystyle{siam}
\bibliography{ref}

\begin{thebibliography}{10}

\bibitem{Alladi}
{\sc K.~Alladi and B.~Gordon}, {\em {Vanishing coefficients in the expansion of products of Rogers-Ramanujan type}}, Contemp. Math., 166 (1994), pp.~129--139.

\bibitem{Andrewsvc}
{\sc G.~E. Andrews and D.~M. Bressoud}, {\em Vanishing coefficients in infinite product expansions}, J. Aust. Math. Soc., 27 (1979), pp.~199--202.

\bibitem{baruah2020some}
{\sc N.~D. Baruah and M.~Kaur}, {\em Some results on vanishing coefficients in infinite product expansions}, Ramanujan J., 53 (2020), pp.~551--568.

\bibitem{C16}
{\sc B.~C. Berndt}, {\em Ramanujan’s notebooks: Part III}, Springer Science \& Business Media, 2012.

\bibitem{cao2011integer}
{\sc Z.~Cao}, {\em Integer matrix exact covering systems and product identities for theta functions}, Int. Math. Res. Notices, 19 (2011), pp.~4471--4514.

\bibitem{keerthana2024generalization}
{\sc Channabasavayya, G.~K. Keerthana, and D.~Ranganatha}, {\em On a generalization of 5-dissections of some infinite $q$-products}, The Ramanujan Journal,  (2024), pp.~1--21.

\bibitem{Ranganatha}
{\sc Channabasavayya and D.~Ranganatha}, {\em On a generalization of vanishing coefficients in two $q$-series expansions}, J. Ramanujan Math. Soc., 39 (2024), pp.~187--–192.

\bibitem{dasappa2024vanishing}
\leavevmode\vrule height 2pt depth -1.6pt width 23pt, {\em Vanishing coefficients in arithmetic progressions of some infinite products with modulo 13, 17, 19 and 29}, The Journal of Analysis,  (2024), pp.~1--14.

\bibitem{hirschhorn2019two}
{\sc M.~D. Hirschhorn}, {\em Two remarkable $q$-series expansions}, Ramanujan J., 49 (2019), pp.~451--463.

\bibitem{vm2}
{\sc M.~Kaur and Vandna}, {\em Results on vanishing coefficients in infinite $q$-series expansions for certain arithmetic progressions mod 7}, Ramanujan J., 58 (2022), pp.~269--289.

\bibitem{liu2024vanishing}
{\sc J.-C. Liu}, {\em On the vanishing coefficients of odd powers of {R}amanujan’s theta functions}, Ramanujan J., 65 (2024), pp.~45--52.

\bibitem{Mc3}
{\sc J.~Mc~Laughlin}, {\em {A generalization of Schr{\"o}ter’s formula}}, Anna. Comb., 23 (2019), pp.~889--906.

\bibitem{mc2021new}
\leavevmode\vrule height 2pt depth -1.6pt width 23pt, {\em New infinite $q$-product expansions with vanishing coefficients}, Ramanujan J., 55 (2021), pp.~733--760.

\bibitem{ramanujan1988lost}
{\sc S.~Ramanujan}, {\em The lost notebook and other unpublished papers}, Narosa Pub. House, 1988.

\bibitem{richmond1978taylor}
{\sc B.~Richmond and G.~Szekeres}, {\em {The Taylor coefficients of certain infinite products}}, Acta Sci. Math., 40 (1978), pp.~347--369.

\bibitem{tang2019vanishing}
{\sc D.~Tang}, {\em Vanishing coefficients in some $q$-series expansions}, Int. J. Number Theory, 15 (2019), pp.~763--773.

\bibitem{2tang2023vanishing}
\leavevmode\vrule height 2pt depth -1.6pt width 23pt, {\em Vanishing coefficients in three families of products of theta functions}, Revista de la Real Academia de Ciencias Exactas, F{\'\i}sicas y Naturales. Serie A. Matem{\'a}ticas, 117 (2023), p.~36.

\bibitem{tang2023vanishing}
\leavevmode\vrule height 2pt depth -1.6pt width 23pt, {\em {Vanishing coefficients in three families of products of theta functions-II}}, Results Math., 78 (2023).

\bibitem{vm1}
{\sc Vandna and M.~Kaur}, {\em Vanishing coefficients of $q^{5n+r}$ and $q^{11n+r}$ in certain infinite $q$-product expansions}, Anna. Comb., 26 (2022), pp.~533--557.

\end{thebibliography}

\end{document}